\documentclass[a4paper,12pt]{article}
\usepackage{amssymb,amsmath,mathrsfs,euscript,amsopn}

\numberwithin{equation}{section}

\newtheorem{lemma}{Lemma}[section]

\newtheorem{theorem}{Theorem}[section]
\newtheorem{proposition}{Proposition}[section]
\newtheorem{corollary}{Corollary}[section]
\newtheorem{definition}{Definition}[section]
\newtheorem{remark}{Remark}[section]

\newenvironment{proof}{\smallskip\noindent{\bf Proof.}\rm}{\hspace*{\fill} $\Box$\medskip}

\newenvironment{prooflemma24}{\smallskip\noindent{\bf Proof of Lemma~\ref{le.8a}.}\rm}{\hspace*{\fill} $\Box$\medskip}
\newenvironment{proofpropa1}{\smallskip\noindent{\bf Proof of Proposition~\ref{appProp}.}\rm}{\hspace*{\fill} $\Box$\medskip}

\title{On the accelerants of non-self-adjoint Dirac operators}

\author{
Ya.~V.~Mykytyuk, D.~V.~Puyda\thanks{\emph{Email addresses:} yamykytyuk@yahoo.com (Ya.~V.~Mykytyuk), dpuyda@gmail.com (D.~V.~Puyda)}\\
\emph{Ivan Franko National University of Lviv}\\
\emph{1 Universytetska str., Lviv, 79000, Ukraine}
}

\date{}

\begin{document}

\maketitle

\begin{abstract}
We prove that there is a homeomorphism between the space of accelerants and the space of potentials of non-self-adjoint Dirac operators on a finite interval.
\end{abstract}

\section{Introduction and main results}

The theory of accelerants was founded by M.~G.~Krein in the middle of the past century. The origins of this theory go back to Krein's short papers \cite{krein1,krein2,krein3,krein4}, where he showed that the resolvent kernels of some integral equations generate solutions of some 2nd order differential equations and systems of 1st order differential equations. Thereby, Krein established a fundamental connection between a special class of functions called the accelerants and Sturm--Liouville and Dirac operators. A detailed presentation of some of these his results can be found in the book \cite{kreinvolterra}.
Krein's ideas in the theory of accelerants were continued and further developed in many papers.

Accelerants play a particular role in the theory of continuous analogues of polynomials orthogonal on the unit circle (see \cite{krein1}). In this context, it is worth mentioning, e.g., remarkable lecture notes \cite{denisov} by S.~A.~Denisov, where the detailed exposition of many aspects of the theory can be found and some new results are obtained.

Let $\mathcal M_r$ denote the Banach algebra of all $r\times r$ matrices with complex entries which we identify with the Banach algebra of linear operators in $\mathbb C^r$ endowed with the standard norm.

\begin{definition}\label{accDef}
We say that a function $h\in L_1((-1,1),\mathcal M_r)$ is an accelerant if for each $\alpha\in(0,1]$ the integral equation
\begin{equation}\label{accEq}
f(x)+\int_0^\alpha h(x-t)f(t)\,{\mathrm d} t=0,\qquad x\in(0,1),
\end{equation}
has only zero solution in $L_2((0,1),\mathbb C^r)$.
\end{definition}

Note that Definition~\ref{accDef} differs from the one originally introduced by Krein in that we do not require any of the conditions $h(x)=h(-x)$ or $h(x)=h(-x)^*$, $x\in(-1,1)$.
Note also that if $h$ is an accelerant, then such is also $h^\sharp$, where $h^\sharp(x):=h(-x)$, $x\in(-1,1)$ (see Remark~\ref{rem3} below).

We denote by $\mathfrak H_{p,r}$ the set of accelerants belonging to
$L_p((-1,1),\mathcal M_r)$, $p\in[1,\infty)$, and endow $\mathfrak H_{p,r}$ with the metric of the latter.
It is known (see Proposition~\ref{hCor} below) that for an arbitrary accelerant $h\in\mathfrak H_{p,r}$, the integral equation
\begin{equation}\label{kreinEq}
r(x,t)+h(x-t)+\int_0^x r(x,s)h(s-t)\,{\mathrm d} s=0,\qquad (x,t)\in\overline\Omega_+,
\end{equation}
where $\Omega_+:=\{(x,t)\mid0<t<x<1\}$, has a unique solution $r_h\in L_1(\Omega_+,\mathcal M_r)$. If one sets $r_h(x,t)=0$ for $(x,t)\in[0,1]^2\setminus\Omega_+$, then $r_h\in G_{p,r}^+$ (see definition in Sect.~\ref{Gsect} below). Equation (\ref{kreinEq}) is called the \emph{Krein equation}.

The connection between the accelerants and Dirac systems of differential equations was established by Krein in \cite{krein4}. In the present paper, it is convenient to explain this connection in equivalent form using the solution of equation (\ref{kreinEq}).

So, let $h\in\mathfrak H_{p,r}$. Consider the $r\times r$ matrix-valued functions
$$
\varphi_1(x,\lambda):={\mathrm e}^{{\mathrm i}\lambda x}
\left(I+\int_0^x {\mathrm e}^{-2{\mathrm i}\lambda s}r_h(x,x-s)\,{\mathrm d} s\right),
$$
$$
\varphi_2(x,\lambda):={\mathrm e}^{-{\mathrm i}\lambda x}
\left(I+\int_0^x {\mathrm e}^{2{\mathrm i}\lambda s}r_{h^\sharp}(x,x-s)\,{\mathrm d} s\right),
$$
where $x\in(0,1)$, $\lambda\in\mathbb C$, $I$ is the $r\times r$ identity matrix and $r_{h^\sharp}$ is the solution of (\ref{kreinEq}) with $h^\sharp$ instead of $h$. Then the $2r\times r$ matrix-valued function $\varphi:=(\varphi_1,\,\varphi_2)^\top$ is a solution of the Cauchy problem
$$
J\frac{\mathrm d}{\mathrm d x}\varphi+Q\varphi=\lambda\varphi,\qquad
\varphi(0,\lambda)=\begin{pmatrix}I\\I\end{pmatrix},
$$
with
$$
J:=\frac1{\mathrm i}\begin{pmatrix}I&0\\0&-I\end{pmatrix},\qquad
Q(x)=[\Theta(h)](x):=\begin{pmatrix}0&{\mathrm i} r_h(x,0)\\-{\mathrm i} r_{h^\sharp}(x,0)&0\end{pmatrix},\quad x\in(0,1).
$$
Since both functions $r_h$ and $r_{h^\sharp}$ belong to $G_{p,r}^+$, the function $Q=\Theta(h)$ belongs to the class
$$
\mathfrak Q_p:=\{Q\in L_p((0,1),\mathcal M_{2r})\mid Q(x)J=-JQ(x)\,\,\text{a.e. on}\,\,(0,1)\}.
$$
The mapping $\Theta:\mathfrak H_{p,r}\to\mathfrak Q_p$ will be called the \emph{Krein mapping}.

The main result of this paper is the following theorem:

\begin{theorem}\label{mainTh}
For an arbitrary $p\in[1,\infty)$, the Krein mapping is a homeomorphism between the metric spaces $\mathfrak H_{p,r}$ and $\mathfrak Q_p$. Moreover, both the Krein mapping and its inverse are locally Lipschitz.
\end{theorem}

In his paper \cite{krein4}, Krein treated symmetric accelerants, i.e. the ones satisfying the condition $h(-t)=h(t)=h(t)^\top$, $t\in(-1,1)$, where $\top$ designates the transposition of matrices. Namely, he proved that there is a one-to-one correspondence between the set of all continuous symmetric accelerants and the set of all continuous symmetric potentials of the Krein systems which are closely related to Dirac operators.

The analogue of Krein's theorem was established for self-adjoint Dirac operators with continuous potentials in \cite{5auth}. Therein, it was shown that there is a one-to-one corres\-pondence between the potentials of such operators and hermitian accelerants (i.e. such that $h(-t)=h(t)^*$) that are continuous outside the origin.

The analogous theorem for Krein systems on semi-axis was proved in \cite[Theorem~5.3]{denisov}.

In \cite[Theorem~1.9]{sturm}, it was proved that the Krein mapping is a homeomorphism between the space of all even accelerants $h$ in $\mathfrak H_{2,r}$ and $L_2((0,1),\mathcal M_r)$. In \cite[Theorem~1.5]{dirac2}, the same result was established about the space of hermitian accelerants $h$ in $\mathfrak H_{p,r}$ and $L_p((0,1),\mathcal M_r)$, $p\in[1,\infty)$.
It thus follows that the potentials of all \emph{self-adjoint} Dirac operators on $[0,1]$ correspond to \emph{hermitian} accelerants $h\in\mathfrak H_{p,r}$. In the present paper, we actually abandon the condition of self-adjointness and show that the potentials of all (not necessarily self-adjoint) Dirac operators on $[0,1]$ correspond to (not necessarily hermitian) accelerants $h\in\mathfrak H_{p,r}$. Since $r_{h^\sharp}(\cdot,0)=[r_h(\cdot,0)]^*$ for all hermitian accelerants, the results of the present paper correlate well with the results of \cite{dirac2}.

\section{Some facts from the theory of factorizations}

\subsection{Some general facts}

Let $\mathcal H$ be a separable infinite dimensional Hilbert space and $\mathcal B:=\mathcal B(\mathcal H)$ be the Banach algebra of all everywhere defined bounded linear operators in $\mathcal H$. We write $\mathcal B_\infty$ and $\mathcal B_0$ for the Banach algebra of all compact operators and for the linear space of all finite dimensional operators from $\mathcal B$, respectively.

We say that a set $\mathfrak P\subset\mathcal B$ of orthoprojectors is a \emph{chain} if for any $P_1,P_2\in\mathfrak P$ it holds either $P_1<P_2$ or $P_2<P_1$. A chain is said to be \emph{closed} if it is a closed subset of $\mathcal B$ in the strong operator topology. A closed chain is said to be \emph{continuous} if for each pair $P_1,P_2\in\mathfrak P$ such that $P_1<P_2$ there is
$P\in\mathfrak P$ such that $P_1<P<P_2$. We say that a closed chain $\mathfrak P$ is \emph{complete} if it is continuous and $0,I\in\mathfrak P$, where $I$ is the identity operator in $\mathcal H$.

Let $\mathfrak P$ be a complete chain in $\mathcal H$. Set
\begin{align*}
{{\mathcal B}_\infty^+} &:=\{B\in {\mathcal B}_\infty \mid \forall P\in\mathfrak P \quad (I-P)BP=0 \},
\\
{{\mathcal B}_\infty^-} &:=\{B\in {\mathcal B}_\infty \mid \forall P\in\mathfrak P \quad  PB(I-P)=0\}.
\end{align*}
It can be easily verified that $\mathcal B_\infty^+$ and $\mathcal B_\infty^-$ are closed Banach subalgebras in $\mathcal B_\infty$ and that $\mathcal B_\infty^+\cap\mathcal B_\infty^-=\{0\}$. Furthermore, the operators from $\mathcal B_\infty^\pm$ are Volterra ones (see~\cite[Ch.~I]{kreinvolterra}).

Denote by $\mathcal P^+$ ($\mathcal P^-$, resp.) the projector in $\widetilde{\mathcal B}_\infty:=\mathcal B_\infty^+\dot+\mathcal B_\infty^-$ onto $\mathcal B_\infty^+$ ($\mathcal B_\infty^-$, resp.) parallel to $\mathcal B_\infty^-$ ($\mathcal B_\infty^+$, resp.). The projectors $\mathcal P^+$ and $\mathcal P^-$ are called the \emph{transformators of triangular truncations} (this term was suggested by I.~C.~Gohberg and M.~G.~Krein for the operators acting from one Banach algebra to another, see~\cite[Ch.~II]{kreinvolterra}).

Denote by $\Sigma$ the set of all Banach algebras $\mathfrak S\subset\mathcal B_\infty$ in which the transformators $\mathcal P^+$ and $\mathcal P^-$ are continuous. For each $\mathfrak S\in\Sigma$ we set
\begin{equation}\label{pm1}
\mathfrak S^\pm:=\mathcal P^\pm\mathfrak S.
\end{equation}
It then follows that both $\mathfrak S^+$ and
$\mathfrak S^-$ are closed subalgebras in $\mathfrak S$ consisting of Volterra operators and that
$\mathfrak S=\mathfrak S^+\dot+\mathfrak S^-$.

Let $\mathfrak S\in\Sigma.$  We say that the operator $I+Q$ with $Q\in\mathcal B_\infty$
($Q\in\mathfrak S$, resp.), \emph{admits a factorization} in $\mathcal B_\infty$ (in $\mathfrak S$, resp.) if
\begin{equation}\label{Qfact}
I+Q=(I+K_-)^{-1}(I+K_+)^{-1}
\end{equation}
with some $K_{\pm}\in\mathcal B^\pm_\infty$ ($K_{\pm}\in\mathfrak S^\pm $, resp.).

Let $\Phi$ ($\Phi_{\mathfrak S}$, resp.) denote the set of all operators $Q\in\mathcal B_\infty$
 ($Q\in\mathfrak S$) for which $I+Q$ admits a factorization in $\mathcal B_\infty$ (in $\mathfrak S $). It is known (see~\cite[Ch.~IV]{kreinvolterra}) that $\Phi$ is contained in the set
\begin{equation*}
\Psi := \{Q \in {\mathcal B}_{\infty}\mid
\forall P\in\mathfrak P \quad  \ker (I+PQP)=\{0\}\}
\end{equation*}
and that for each $Q\in\Phi$ the operators $K_\pm=K_\pm(Q)$ in (\ref{Qfact}) are determined uniquely.
 The following theorem is proved in \cite{MykFact1}:

\begin{theorem}\label{th.1a}
Let $\mathfrak S\in\Sigma$. Then the set $\Phi_{\mathfrak S}$ is open in $\mathfrak S.$  Moreover, the mappings $\Phi_{\mathfrak S}\owns Q\mapsto K_\pm(Q)\in\mathfrak S$ are locally Lipschitz.
\end{theorem}

\begin{remark}\label{re.2a}\rm
A mapping $\varphi$ acting from an open set $\mathcal O$ in a Banach space $X$ to a Banach space $Y$ is said to be \emph{locally Lipschitz} if for each $x_0\in\mathcal O$ there are a neighbourhood $\mathcal U\subset\mathcal O$ of $x_0$ and $c>0$ such that $\|\varphi(x_1)-\varphi(x_2)\|_Y\le c\|x_1-x_2\|_X$ for all $x_1,x_2\in\mathcal U$.
\end{remark}

Set
$$    \Sigma_f:=\{\mathfrak S\in\Sigma \mid \Phi_{\mathfrak S}=
\Psi\cap{\mathfrak S}\}, \quad \Sigma_f^0:=\{\mathfrak S\in\Sigma_f \mid {\mathfrak S}\cap {\mathcal B}_0  \,\,\text{is dense everywhere in} \,\,{\mathcal B}_{\infty}\}.
$$
Note that as follows from the well known results in the theory of factorizations (see~\cite{kreinvolterra}) the Neumann--Schatten ideals $\mathcal B_p$, $1<p<\infty$, belong to the class $\Sigma_f^0$.

The next two theorems follow from the results of \cite{MykFact1,MykFact2}:

\begin{theorem}\label{th.3a}
Let $\mathfrak S\in\Sigma$ and $\mathfrak S_1\in\Sigma_f^0$ be a two-sided ideal in $\mathfrak S$. If $\mathfrak S_1$ is dense everywhere in $\mathfrak S$, then $\mathfrak S\in\Sigma_f^0$.
\end{theorem}

\begin{theorem}\label{th.4a}
Let  $Q\in\Phi$ and $Q_1\in\mathcal B_0.$  Then the set $\{ \lambda\in\mathbb C \mid (Q+\lambda Q_1)\in \Psi\}$ is open and  dense everywhere in $\mathbb C$.
\end{theorem}

\begin{corollary}\label{cor.4a}
Let $\mathfrak S\in\Sigma_f^0$. Then the set $\Phi_{\mathfrak S}$ is dense everywhere in $\mathfrak S$.
\end{corollary}

\subsection{Algebras $\mathscr G_{p,n}$}\label{Gsect}

For an arbitrary $p\in[1,\infty)$ and $n\in\mathbb N$, we denote by $G_{p,n}$ the set of all measurable functions $K:[0,1]^2\to\mathcal M_n$ such that for all $x,t\in[0,1]$ the functions $K(x,\cdot)$ and $K(\cdot,t)$ belong to $L_p((0,1),\mathcal M_n)$ and, moreover, the mappings
$$
[0,1]\ni x\mapsto K(x,\cdot)\in L_p((0,1),\mathcal M_n), \qquad
[0,1]\ni t\mapsto K(\cdot,t)\in L_p((0,1),\mathcal M_n)
$$
are continuous. The set $G_{p,n}$ becomes a Banach space upon introducing the norm
\begin{equation}\label{Gnorm}
\|K\|_{G_{p,n}}=\max\left\{
\max\limits_{x\in[0,1]}\|K(x,\cdot)\|_{L_p},\
\max\limits_{t\in[0,1]}\|K(\cdot,t)\|_{L_p}\right\}.
\end{equation}
We denote by $\mathscr G_{p,n}$ the set of all integral operators in $\mathcal H:=L_2((0,1),\mathbb C^n)$ with kernels $K\in G_{p,n}$ and endow $\mathscr G_{p,n}$ with the norm
$$
\|\mathscr K\|_{\mathscr G_{p,n}}:=\|K\|_{G_{p,n}},\qquad
\mathscr K\in\mathscr G_{p,n}.
$$

Note that there are continuous embeddings $\mathscr G_{p,n}\subset \mathscr G_{1,n} \subset \mathcal B(\mathcal H)$ and that for each $\mathscr K\in \mathscr G_{p,n}$ and $\mathscr R\in \mathscr G_{1,n}$ it holds
$$
\|\mathscr K\|_{\mathscr G_{1,n}}\le\|\mathscr K\|_{\mathscr G_{p,n}}, \qquad
\|\mathscr R\|_{\mathcal B}\le\|\mathscr R\|_{\mathscr G_{1,n}}.
$$
Furthermore, it can be verified that $\mathscr G_{p,n}$ is a Banach algebra.

We set
$$
\quad \Omega_+:=\{(x,t)\mid0<t<x<1\},\qquad \Omega_-:=\{(x,t)\mid0<x<t<1\}
$$
and write $G_{p,n}^\pm$ for the sets of all functions $K\in G_{p,n}$ such that $K(x,t)=0$ a.e. in $\overline\Omega_\mp$. We denote by $\mathscr G_{p,n}^\pm$ the subalgebras in $\mathscr G_{p,n}$ consisting of all operators $\mathscr K\in\mathscr G_{p,n}$ with kernels $K\in G_{p,n}^\pm$.
It is easy to verify that $\mathscr G_{p,n}^\pm$ are closed subalgebras in
$\mathscr G_{p,n}$ and that $\mathscr G_{p,n}=\mathscr G_{p,n}^+\dotplus \mathscr G_{p,n}^-$.

We denote by $\mathscr S^\pm_n$ the operator algebras consisting of all operators $\mathscr K\in\mathscr G_{1,n}^\pm$ with kernels that are continuous in $\overline\Omega_\pm$.
The algebras $\mathscr S_n^+$ and $\mathscr S_n^-$ become Banach algebras upon introducing the norms
$$
\|\mathscr K\|_{\mathscr S_n^\pm}:=  \max_{(x,t)\in\overline\Omega_\pm}\|K(x,t)\|.
$$
We set $\mathscr S_n:=\mathscr S_n^+ \dotplus\mathscr S_n^-$ and endow $\mathscr S_n$ with the norm
$$
   \|\mathscr K\|_{\mathscr S_n}:=  \max\left\{\|\mathscr K_+\|_{\mathscr S_n^+},\,\,
    \|\mathscr K_-\|_{\mathscr S_n^-} \right\}, \qquad \mathscr K=\mathscr K_+ +\mathscr K_-,
    \quad \mathscr K_\pm\in \mathscr S_n^\pm.
$$
It is easy to verify that $\mathscr S_n$ is a Banach algebra. We then denote by $\mathscr S_{n,0}$ a subalgebra in $\mathscr S_n$ consisting of all operators $\mathscr K\in\mathscr S_n$ with kernels that are continuous on $[0,1]^2$.

\begin{lemma}\label{le.4a}
$\mathscr S_n$ is a two sided ideal in $\mathscr G_{p,n}$. Furthermore, $\mathscr S_n$ and $\mathscr S_n^\pm$ are continuously and densely embedded into $\mathscr G_{p,n}$ and $\mathscr G_{p,n}^\pm$, respectively.
\end{lemma}
\begin{proof}
A straightforward verification shows that for each $\mathscr K_0\in \mathscr S_n$ and
$\mathscr K_1\in \mathscr G_{p,n}$, the products $\mathscr K_0 \mathscr K_1$ and
$\mathscr K_1 \mathscr K_0$ belong to $ \mathscr S_{n,0}$ and that
\begin{equation}\label{eq.22y}
 \|\mathscr K_0 \mathscr K_1\|_{\mathscr S_n},\|\mathscr K_1 \mathscr K_0\|_{\mathscr S_n}\le \|\mathscr K_0 \|_{\mathscr S_n}
 \|\mathscr K_1\|_{\mathscr G_{p,n}}.
\end{equation}
Therefore, $\mathscr S_n$ is a two sided ideal in $\mathscr G_{p,n}$. Since
\begin{equation}\label{eq.23y}
  \|\mathscr K\|_{\mathscr G_{p,n}}\le \|\mathscr K\|_{\mathscr S_n}, \qquad\mathscr K\in\mathscr S_n,
\end{equation}
one also has that $\mathscr S_n$ and $\mathscr S_n^\pm$ are continuously embedded into $\mathscr G_{p,n}$ and $\mathscr G_{p,n}^\pm$, respectively.

It was proved in \cite{MykFact2} that $\mathscr S_{1,0}$ is dense everywhere in
$\mathscr G_{1,1}$. By a straightforward modification of that proof it can be shown
that $\mathscr S_{n,0}$ is dense everywhere in $\mathscr G_{p,n}$. It then follows that
$\mathscr S^+_n$ and $\mathscr S^-_n$ are dense everywhere in $\mathscr G^+_{p,n}$ and $\mathscr G^-_{p,n}$, respectively.
\end{proof}

\noindent In particular, it follows from Lemma~\ref{le.4a} that $\mathscr G_{p,n}\subset\mathcal B_\infty(\mathcal H)$.

\begin{lemma}\label{le.5a}
Let $\mathscr K\in\mathscr G_{p,n}^+\cup\mathscr G_{p,n}^-$ and $\rho(\mathscr K)$ be the spectral radius of $\mathscr K$ (see \cite[Chap.~10]{rudin}). Then $\rho(\mathscr K)=0$.
\end{lemma}

\begin{proof}
Since the mapping $\mathscr K\mapsto \mathscr K^* $ maps $\mathscr G_{p,n}^-$ onto $\mathscr G_{p,n}^+$ isometrically, it suffices to prove that $\rho(\mathscr K)=0$ for each $\mathscr K\in\mathscr G_{p,n}^+$.

For this purpose, note that for an arbitrary sequence $(\mathscr K_j)_{j=1}^m$ in $\mathscr S_n^+$ it holds
\begin{equation}\label{eq.22z}
\|\mathscr K_1 \cdots \mathscr K_m\|_{\mathscr S_n}   \le \frac1{m!} \prod_{j=1}^m
\|\mathscr K_j\|_{\mathscr S_n}.
\end{equation}
Let $\mathscr K\in\mathscr G_{p,n}^+$ and $\delta\in(0,1)$.
In view of Lemma~\ref{le.4a}, the operator $\mathscr K$ can be written in the form $\mathscr K=\mathscr K_0+\mathscr K_1$ with some
$\mathscr K_0\in\mathscr S^+_n$ and  $\mathscr K_1\in\mathscr G_{p,n}^+$ such that $\|\mathscr K_1\|_{\mathscr G_{p,n}} \le \delta$.
It then holds
$\mathscr K^s =\sum\limits_{\sigma\in  U_s} \mathscr K_{\sigma(1)} \cdots \mathscr  K_{\sigma(s)}$, where the sum is taken over the set
$U_s$ of all functions $\sigma:\{1,\dots,s\}\to \{0,1\}$, and thus one has
$$
\|\mathscr K^s\|_{\mathscr G_{p,n}} \le 2^s \max_{\sigma} \|\mathscr K_{\sigma(1)} \cdots \mathscr  K_{\sigma(s)}\|_{\mathscr G_{p,n}}.
$$
Let $\sigma\in U_n$ and $m:={\operatorname{card}}(\sigma^{-1}(0))>0$. It then follows from Lemma~\ref{le.4a} and from the estimates (\ref{eq.22y}) and (\ref{eq.22z}) that the operator
$\mathscr L=\mathscr K_{\sigma(1)} \cdots \mathscr  K_{\sigma(s)}$
belongs to $\mathscr S^+_n$ and, furthermore,
$$
\|\mathscr L\|_{\mathscr S_n}\le
\frac1{m!}{\|\mathscr K_1\|^{s-m}_{\mathscr G_{p,n}} \|\mathscr K_0\|^{m}_{\mathscr S_n}}
\le\frac1{m!}{\delta^{s-m}  \|\mathscr K_0\|_{\mathscr S_n}^m }.
$$
Taking into account (\ref{eq.23y}), we then obtain that
\begin{equation*}\label{eq.24y}
  \|\mathscr L\|_{\mathscr G_{p,n}}\le
 \delta^s \frac{(\delta^{-1} \|\mathscr K_0\|_{\mathscr S_n})^m }{m!}\le
  \delta^s \exp{(\delta^{-1} \|\mathscr K_0\|_{\mathscr S_n})} .
\end{equation*}
Evidently, the latter inequality holds true also for $m=0$.
Therefore, one has
$$
\|\mathscr K^s\|_{\mathscr G_{p,n}}\le (2\delta)^s \exp{(\delta^{-1} \|\mathscr K_0\|_{\mathscr S_n})}
$$
and thus the spectral radius $\rho(\mathscr K)$ of the operator $\mathscr K\in\mathscr G_{p,n}$
does not exceed $2\delta$. Since $\delta$ was arbitrary, one then has that $\rho(\mathscr K)=0$.
\end{proof}

It follows from Lemma~\ref{le.5a} that the mapping
$
\mathscr K \mapsto \gamma(\mathscr K):=(I+\mathscr K)^{-1}-I
$
maps both $\mathscr G_{p,n}^+$ and $\mathscr G_{p,n}^-$ into themselves. Actually even more holds true:

\begin{lemma}\label{le.7a}
The mappings
$\mathscr G_{p,n}^+ \ni\mathscr K \mapsto \gamma(\mathscr K)\in \mathscr G_{p,n}^+ $ and $\mathscr G_{p,n}^- \ni\mathscr K \mapsto \gamma(\mathscr K)\in \mathscr G_{p,n}^- $ are homeomorphic and locally Lipschitz.
\end{lemma}

\noindent Lemma~\ref{le.7a} follows from the next general result:

\begin{proposition}\label{pr.6a}
Let $\mathcal A$ be a Banach algebra with the identity $e$ and $\mathcal A_0$ be its closed subalgebra such that $\rho(a)=0$ for each $a\in\mathcal A_0$. Then the mapping
$
\mathcal A_0 \ni a \mapsto \gamma(a):=[(e+a)^{-1}-e]\in\mathcal A_0
$
is homeomorphic and locally Lipschitz.
\end{proposition}

\begin{proof}
Since $(e+a_1)^{-1} -(e+a_2)^{-1}=(e+a_1)^{-1}(a_2-a_1)(e+a_2)^{-1}$, it follows that
\begin{equation}\label{eq.26y}
\gamma(a_1) -\gamma(a_2)=(e +\gamma(a_1))(a_2-a_1)(e+\gamma(a_2)), \qquad  a_1, a_2\in\mathcal A.
\end{equation}
Let $a\in\mathcal A_0$. It then follows from the assumptions of the lemma that there is $m\in\mathbb N$ such that $\|a^m\|_{\mathcal A}\le 1/4$. Set $C:=\sum_{k=0}^{m-1} \|a^k\|$.
Since multiplication in $\mathcal A_0$ is continuous, it then follows that there is a neighbourhood
$\mathcal U\subset\mathcal A_0$ of $a$ such that
\begin{equation*}
\|b^m\|_{\mathcal A}\le 1/2, \qquad  \sum_{k=0}^{m-1} \|b^k\|_{\mathcal A} \le 2C,
 \qquad b\in\mathcal U.
\end{equation*}
Since
$$  \gamma(b)=\sum\limits_{s=1}^\infty (-b)^s =\sum\limits_{k=0}^{m-1}(-b)^k \sum\limits_{s=1}^\infty (-b)^{ms},
$$
it follows that
$  \|\gamma(b)\|_{\mathcal A}\le  2C $ for all $b\in\mathcal U$.
Taking into account (\ref{eq.26y}) we then obtain that
\begin{equation*}
\|\gamma(a_1) -\gamma(a_2)\|_{\mathcal A}\le (1+2C)^2\|a_1-a_2\|_{\mathcal A}, \qquad a_1,a_2\in \mathcal U.
\end{equation*}
Therefore  the mapping $\gamma$ is locally Lipschitz. Since $\gamma(\gamma(a))\equiv a$, it also follows that $\gamma$ is homeomorphic.
\end{proof}

\subsection{Factorization of operators in $\mathscr G_{p,n}$}

Let $\mathcal H:=L_2((0,1),\mathbb C^n)$. We consider the transformators $\mathcal P^\pm$ in $\mathcal B_\infty(\mathcal H)$ generated by a complete chain of orthoprojectors $P_\alpha:\mathcal H\to\mathcal H$, $\alpha\in [0,1]$, given by the formula
$$
P_\alpha f:=\chi_{[0,\alpha]} f,  \qquad f\in\mathcal H,
$$
where $\chi_{[0,\alpha]} $ is the characteristic function of the interval $[0,\alpha]$.

\begin{lemma}\label{le.8a}
The transformators $\mathcal P^\pm$ are continuous in $\mathscr G_{p,n}$ and $\mathscr G_{p,n}^\pm=\mathcal P^\mp\mathscr G_{p,n}$.
\end{lemma}

\begin{remark}\rm
The change of sign in the above formula arises due to discrepancy between definition (\ref{pm1}) and the definition of algebras $\mathscr G_{p,n}^\pm$ given at the beginning of Sect.~\ref{Gsect}. However, the authors prefer to accept this inconvenience in order to follow both the standard notations in \cite{kreinvolterra,MykFact1,MykFact2} and the ones used in \cite{dirac1,sturm,dirac2}. This causes also sign differences between formula (\ref{eq.27y}) below and formula (\ref{Qfact}).
\end{remark}

\begin{prooflemma24}
In the scalar case $n=1$, continuity of  $\mathcal P_1^\pm:=\mathcal P^\pm$ in $\mathscr G_{p,1}$ follows from the results of \cite{MykFact2}. Note that $\mathscr G_{p,n}$ can be considered as a tensor product of the algebras $\mathscr G_{p,1}$ and $\mathcal M_n$ and that $\mathcal P^\pm $ can be considered as tensor products of the operators $\mathcal P_1^\pm$ and $I_{\mathcal M_n}$. Therefore, we obtain that the transformators $\mathcal P^\pm$ act continuously in $ \mathscr G_{p,n}$.
Verification of the equalities $\mathscr G_{p,n}^\pm=\mathcal P^\mp\mathscr G_{p,n}$ is straightforward.
\end{prooflemma24}

\begin{remark}\rm\label{remA}
It follows from Lemma~\ref{le.8a} that $\mathscr G_{p,n}$ belongs to the class $\Sigma$.
The algebra $\mathscr S_n$ belongs to the class $\Sigma_f$
(see~\cite[Ch.~IV]{kreinvolterra}). Since $\mathscr S_{n}\cap \mathcal B_0(\mathcal H)$ is dense everywhere in $\mathcal B_\infty(\mathcal H)$, it follows that
$\mathscr S_{n}$ belongs to the class $\Sigma_f^0$. Therefore, taking into account
Theorem~\ref{th.3a} and Lemma~\ref{le.4a}, we obtain that $\mathscr G_{p,n}$ also belongs to $\Sigma_f^0$.  \end{remark}

We denote by $\widetilde {\mathscr G}_{p,n}$ the set of all operators $\mathscr F\in \mathscr G_{p,n}$ such that $\mathscr I+\mathscr F$ admits a factorization in $\mathscr G_{p,n}$, i.e.
$\mathscr F\in \widetilde {\mathscr G}_{p,n}$ if and only if
there exist $\mathscr L_\pm\in \mathscr G_{p,n}^\pm$ such that
\begin{equation}\label{eq.27y}
\mathscr I+\mathscr F=(\mathscr I+\mathscr L_+)^{-1}(\mathscr I+\mathscr L_-)^{-1}.
\end{equation}
In view of Theorem~\ref{th.1a} and Corollary~\ref{cor.4a} we then arrive at the following statements:

\begin{theorem}\label{th9b}
\begin{itemize}
  \item[(i)]  The set $\widetilde {\mathscr G}_{p,n}$ is open and dense everywhere in $\mathscr G_{p,n}$.
  \item[(ii)] If $\mathscr F\in \widetilde {\mathscr G}_{p,n}$, then $\mathscr L_\pm$ in (\ref{eq.27y}) are determined uniquely and $\mathscr L_\pm=K_\mp(\mathscr F)$.
  \item[(iii)] The mapping $\widetilde {\mathscr G}_{p,n}\ni \mathscr F\mapsto K_{\pm}(\mathscr F)\in \mathscr G_{p,n} $ is locally Lipschitz.
\end{itemize}
\end{theorem}

\begin{theorem}\label{th9a}
Let $\mathscr F\in {\mathscr G}_{p,n}$ and $F\in G_{p,n}$ be a kernel of $\mathscr F$.
Then the following statements are equivalent:
\begin{itemize}
\item[(i)] $\mathscr F\in \widetilde {\mathscr G}_{p,n}$;
\item[(ii)] for each $\alpha\in [0,1]$, the integral equation
\begin{equation}\label{accEq1}
f(x)+\int_0^\alpha F(x,t)f(t)\,{\mathrm d} t=0,\qquad x\in(0,1),
\end{equation}
has only zero solution in $\mathcal H$;
\item[(iii)] the integral equation
\begin{equation}\label{GLMFact}
X(x,t)+F(x,t)+\int_0^xX(x,s)F(s,t)\,{\mathrm d} s=0,\qquad (x,t)\in\overline\Omega_+,
\end{equation}
is solvable in $G^+_{p,n}$.
\end{itemize}
\end{theorem}

\begin{remark}\label{rem1}\rm
Equation (\ref{GLMFact}) always has at most one solution. If $X\in G^+_{p,n}$ is a solution of (\ref{GLMFact}), then $X$ coincides with the kernel of the operator $\mathscr L_+=K_-(\mathscr F)\in\mathscr G_{p,n}^+$.
\end{remark}

\begin{remark}\label{rem3}\rm
Note that for each $\mathscr F\in\mathscr G_{p,n}$ and $P\in\mathfrak P$, the operators $\mathscr I+\mathscr F P$ and $\mathscr I+P\mathscr F P$ are invertible or not simultaneously. Therefore, it follows that equation (\ref{accEq1}) has a non-zero solution in $\mathcal H$ if and only if it has a non-zero solution in $L_2((0,\alpha),\mathbb C^n)$. For the same reason, we have that the functions $h$ and $h^\sharp$ from $L_p((-1,1),\mathcal M_r)$ belong to $\mathfrak H_{p,r}$ or not simultaneously.
\end{remark}

\section{Proof of Theorem~\ref{mainTh}}

The aim of this Section is to prove Theorem~\ref{mainTh}, which is the main result of this paper. Firstly, we shall use the results of the previous section to prove that the Krein mapping is locally Lipschitz. Next, we shall construct a locally Lipschitz mapping $\Upsilon:\mathfrak Q_p\to\mathfrak H_{p,r}$ and show that $\Upsilon=\Theta^{-1}$.

\subsection{The Krein mapping}

Here we shall prove that the Krein mapping is locally Lipschitz. We start with several auxiliary statements which will be useful in subsequent expositions. The first one is a corollary of Theorems~\ref{th9b} and \ref{th9a}:

\begin{proposition}\label{hCor}
Let $h\in L_p((-1,1),\mathcal M_n)$. Consider the operator $\mathscr H\in\mathscr G_{p,n}$ acting by the formula
\begin{equation}\label{H1}
(\mathscr Hf)(x)=\int_0^1 h(x-t)f(t)\,{\mathrm d} t,\qquad f\in L_2((0,1),\mathbb C^n).
\end{equation}
Then the following statements are equivalent:
\begin{itemize}
\item[(i)]$\mathscr H\in\widetilde{\mathscr G}_{p,n}$;
\item[(ii)]$h$ is an accelerant, i.e. $h\in\mathfrak H_{p,n}$;
\item[(iii)]the Krein equation
\begin{equation}\label{kreinEq1}
r(x,t)+h(x-t)+\int_0^x r(x,s)h(s-t)\,{\mathrm d} s=0,\qquad (x,t)\in\overline\Omega_+,
\end{equation}
has a unique solution $r_h\in G_{p,n}^+$.
\end{itemize}
Moreover, the mapping $\mathfrak H_{p,n}\owns h\mapsto r_h\in G_{p,n}^+$ is locally Lipschitz and the set $\mathfrak H_{p,n}$ is open in $L_p((-1,1),\mathbb C^n)$.
\end{proposition}

\begin{proposition}\label{hCor1}
The set  $\mathfrak H_{p,r}$  is dense everywhere in $ L_p((-1,1),\mathcal M_r)$.
\end{proposition}

\begin{proof}
Let  $f\in L_p((-1,1),\mathcal M_r)$. Then $f$ can be written in the form $f=h+h_1$, where $h,h_1\in L_p((-1,1),\mathcal M_r)$, $\|h\|_{L_p}<1$ and $h_1$ is a trigonometric polynomial. Denote by $\mathscr H $ and $\mathscr H_1$ the operators constructed by formula (\ref{H1}) from functions $h$ and $h_1$, respectively. It is easily seen that the operator $\mathscr H_1$ is finite dimensional and that the norm of the operator $\mathscr H$ is less than $1$. Therefore, one has $\mathscr H\in\widetilde{\mathscr G}_{p,r}$. By virtue of Theorem~\ref{th.4a}, the set $\Lambda:=\{ \lambda\in\mathbb C \mid (\mathscr H+\lambda \mathscr H_1)\in \Psi\}$ is open and  dense everywhere in $\mathbb C$. Since the algebra ${\mathscr G}_{p,r}$ belongs to the class $\Sigma_f^0$ (see Remark~\ref{remA}), it follows that $\Lambda=\{ \lambda\in\mathbb C \mid (\mathscr H+\lambda \mathscr H_1)\in \widetilde{\mathscr G}_{p,r}\}$. In view of Proposition~\ref{hCor}, this means that $\Lambda=\{ \lambda\in\mathbb C \mid (h+\lambda h_1)\in \mathfrak H_{p,r}\}$. Therefore, $f$ is a limit point of the set $\mathfrak H_{p,r}$.
\end{proof}

\begin{lemma}\label{FhLemma}
Let $h\in L_p((-1,1),\mathcal M_r)$,
\begin{equation}\label{Fh1}
F^h(x,t):=\frac12\begin{pmatrix}h\left(\frac{x-t}2\right)&h\left(\frac{x+t}2\right)\\
h\left(-\frac{x+t}2\right)&h\left(-\frac{x-t}2\right)\end{pmatrix},\qquad x,t\in(0,1),
\end{equation}
and $\mathscr F^h\in\mathscr G_{p,2r}$ be the integral operator with kernel $F^h$. Then
$h\in\mathfrak H_{p,r}\Longleftrightarrow\mathscr F^h\in\widetilde{\mathscr G}_{p,2r}$.
\end{lemma}

\begin{proof}
In view of Proposition~\ref{hCor}, the lemma will be proved if we show that
\begin{equation}\label{equiv}
\mathscr H\in\widetilde{\mathscr G}_{p,r}\quad\Longleftrightarrow\quad
\mathscr F^h\in\widetilde{\mathscr G}_{p,2r}.
\end{equation}
For this purpose, recall (see Theorem~\ref{th9a} and Remark~\ref{rem3}) that $\mathscr H\in\widetilde{\mathscr G}_{p,r}$ if and only if the equation
\begin{equation}\label{Heq}
f(x)+\int_0^\alpha h(x-t)f(t)\,{\mathrm d} t=0,\qquad x\in(0,\alpha),
\end{equation}
has only zero solution in $L_2((0,\alpha),\mathbb C^r)$. Similarly, one has $\mathscr F^h\in\widetilde{\mathscr G}_{p,2r}$ if and only if the equation
\begin{equation}\label{Feq}
g(x)+\int_0^\alpha F^h(x,t)g(t)\,{\mathrm d} t=0,\qquad x\in(0,\alpha),
\end{equation}
has only zero solution in $L_2((0,\alpha),\mathbb C^{2r})$. Now observe that if $f\in L_2((0,\alpha),\mathbb C^r)$ solves (\ref{Heq}), then
$$
g(x)=\begin{pmatrix}f\left(\frac{\alpha+x}2\right)\\f\left(\frac{\alpha-x}2\right)\end{pmatrix}
$$
solves (\ref{Feq}) and that if $g=(g_1,\,g_2)^\top$ with $g_1,g_2\in L_2((0,\alpha),\mathbb C^r)$ solves (\ref{Feq}), then
$$
f(x)=\begin{cases}g_2(\alpha-2x),&x\in\left(0,\tfrac\alpha2\right),\\
g_1(2x-\alpha),&x\in\left(\tfrac\alpha2,\alpha\right),\end{cases}
$$
solves (\ref{Heq}). Therefore, equations (\ref{Heq}) and (\ref{Feq}) have non-zero solutions simultaneously which proves the equivalence (\ref{equiv}).
\end{proof}

\begin{remark}\label{rem2}\rm
Let $h\in\mathfrak H_{p,r}$. Recall that $h^\sharp(x):=h(-x)$ and set
\begin{equation}\label{H}
H(x):=\begin{pmatrix}h(x)&0\\0&h^\sharp(x)\end{pmatrix},\qquad x\in(-1,1).
\end{equation}
It is then easily verified that for the function
\begin{equation}\label{H2}
R_H(x,t):=\begin{pmatrix}r_h(x,t)&0\\0&r_{h^\sharp}(x,t)\end{pmatrix},\qquad (x,t)\in\overline\Omega_+,
\end{equation}
it holds
\begin{equation}\label{RH}
R_H(x,t)+H(x-t)+\int_0^x R_H(x,s)H(s-t)\,{\mathrm d} s=0,\qquad (x,t)\in\overline\Omega_+,
\end{equation}
and that for
\begin{equation}\label{LH}
L_h(x,t):=\frac12\left\{R_H\left(x,\frac{x+t}2\right)
+R_H\left(x,\frac{x-t}2\right)B\right\},\qquad B:=\begin{pmatrix}0&I\\I&0\end{pmatrix},
\end{equation}
one has
\begin{equation}\label{GLM}
F^h(x,t)+L_h(x,t)+\int_0^x L_h(s,t)F^h(x,s)\,{\mathrm d} s=0,\qquad (x,t)\in\overline\Omega_+.
\end{equation}
\end{remark}

Now we are ready to prove the following lemma which is the main purpose of this subsection:

\begin{lemma}\label{dirLip}
For an arbitrary $p\in[1,\infty)$, the Krein mapping $\Theta:\mathfrak H_{p,r}\mapsto\mathfrak Q_p$ is locally Lipschitz.
\end{lemma}

\begin{proof} Let $h\in\mathfrak H_{p,r}$ and $H$ be as in (\ref{H}).
It then follows from Proposition~\ref{hCor} and from (\ref{RH}) that $H\in\mathfrak H_{p,2r}$ and that the mapping $\mathfrak H_{p,2r}\owns H\mapsto R_H\in G_{p,2r}^+$ is locally Lipschitz.
Note that
\begin{equation}\label{ThRH}
[\Theta(h)](x)=R_H(x,0)BJ,\qquad x\in(0,1),
\end{equation}
where $B$ is from (\ref{LH}). In view of formulas (\ref{H}) and (\ref{ThRH}), it is then easily seen that $\Theta$ is locally Lipschitz.
\end{proof}

\subsection{Construction of the mapping $\Upsilon$}

We now construct the mapping $\Upsilon:\mathfrak Q_p\to\mathfrak H_{p,r}$ that will appear to be the inverse of the Krein mapping.

Let $Q\in\mathfrak Q_p$. For each $\lambda\in\mathbb C$, we denote by $\varphi_Q(x,\lambda)$, $x\in[0,1]$, a $2r\times r$ matrix-valued solution of the Cauchy problem
\begin{equation}\label{phiCauchyProbl}
J\frac{\mathrm d}{\mathrm d x}\varphi+Q\varphi=\lambda\varphi,\qquad \varphi(0,\lambda)=\begin{pmatrix}I\\I\end{pmatrix}.
\end{equation}

\begin{lemma}\label{Kprop}
For each $Q\in\mathfrak Q_p$, there is a unique function $K_Q\in G_{p,2r}^+$ such that for all $x\in[0,1]$ and $\lambda\in\mathbb C$ it holds
\begin{equation}\label{phiTrasfOp}
\varphi_Q(x,\lambda)=\varphi_0(x,\lambda)+\int_0^x K_Q(x,s)\varphi_0(s,\lambda)\ {\mathrm d} s,
\end{equation}
where
$\varphi_0(x,\lambda)$ is a solution of the Cauchy problem (\ref{phiCauchyProbl}) in the free case $Q=0$.
Moreover, the mapping $\mathfrak Q_p\owns Q\mapsto K_Q\in G_{p,2r}^+$ is locally Lipschitz.
\end{lemma}

\begin{proof}
Denote by $Y_Q(\cdot,\lambda)\in W_2^1((0,1),\mathcal M_{2r})$ a $2r\times2r$ matrix-valued solution of the Cauchy problem
$$
J\frac{\mathrm d}{\mathrm d x} Y+QY=\lambda Y,\qquad Y(0,\lambda)=I_{2r}.
$$
It then follows from \cite[Theorem~2.1]{cauchy} that there exist unique functions $P^\pm:=P_Q^\pm$ from $G_{p,2r}^+$ such that for all $x\in[0,1]$ and $\lambda\in\mathbb C$ it holds
\begin{equation}\label{Y}
Y_Q(x,\lambda)={\mathrm e}^{-\lambda xJ}+\int_0^x P^+(x,t){\mathrm e}^{-\lambda(x-2t)J}\,{\mathrm d} t
+\int_0^x P^-(x,t){\mathrm e}^{\lambda(x-2t)J}\,{\mathrm d} t.
\end{equation}
Since $\varphi_Q(x,\lambda)=Y_Q(x,\lambda)a$, where $a:=(I,\,I)^\top$, straightforward manipulations lead us to formula (\ref{phiTrasfOp}) with
\begin{equation}\label{K}
K_Q(x,t)=\tfrac12\left\{P^+\left(x,\tfrac{x-t}2\right)+P^+\left(x,\tfrac{x+t}2\right)B
+P^-\left(x,\tfrac{x-t}2\right)B+P^-\left(x,\tfrac{x+t}2\right)\right\},
\end{equation}
where $B$ is from (\ref{LH}).

Let us prove that the mapping $\mathfrak Q_p\owns Q\mapsto K_Q\in G_{p,2r}^+$ is locally Lipschitz.
It follows from the proof of Theorem~2.8 in \cite{cauchy} that with $\widetilde P_Q(x,t):=P_Q^+\left(x,\tfrac{x-t}2\right)$ it holds
\begin{align}
&\|\widetilde P_{Q_1}(x,\cdot)-\widetilde P_{Q_2}(x,\cdot)\|_{L_p}
\le(1+2\varepsilon){\mathrm e}^{2\varepsilon}\|Q_1-Q_2\|_{L_p},
\\
&\|\widetilde P_{Q_1}(\cdot,t)-\widetilde P_{Q_2}(\cdot,t)\|_{L_p}
\le C\|Q_1-Q_2\|_{L_p},\quad C:=2\varepsilon{\mathrm e}^\varepsilon+2\varepsilon(1+2\varepsilon){\mathrm e}^{2\varepsilon},
\end{align}
for every $Q_1,Q_2\in\mathfrak Q_p$ such that $\|Q_1\|,\|Q_2\|<\varepsilon$ and that
the same estimates hold true also with $\widetilde P_Q(x,t):=P_Q^+\left(x,\tfrac{x+t}2\right)$, $\widetilde P_Q(x,t):=P_Q^-\left(x,\tfrac{x+t}2\right)$ and $\widetilde P_Q(x,t):=P_Q^-\left(x,\tfrac{x-t}2\right)$. In view of (\ref{K}), we then obtain that the mapping $\mathfrak Q_p\owns Q\mapsto K_Q\in G_{p,2r}^+$ is locally Lipschitz.
\end{proof}

Denote by $\mathscr K_Q\in\mathscr G_{p,2r}$ the integral operator with kernel $K_Q$ and let $\mathscr I$ stand for the identity operator in $\mathbb H$. Since $\mathscr K_Q$ is a Volterra operator, the operator $\mathscr I+\mathscr K_Q$ is invertible in $\mathbb H$. Set
\begin{align}\label{L}
\mathscr L_Q&:=(\mathscr I+\mathscr K_Q)^{-1}-\mathscr I,
\\\label{F}
\mathscr F_Q&:=(\mathscr I+\mathscr K_Q)^{-1}(\mathscr I+\mathscr K_{Q^*}^*)^{-1}-\mathscr I
\end{align}
and denote by $L_Q$ and $F_Q$ the kernels of the integral operators $\mathscr L_Q$ and $\mathscr F_Q$, respectively.

\begin{theorem}\label{UpsTh}
Let $Q\in\mathfrak Q_p$ and $F:=F_Q$. Then there is a unique $h=\Upsilon(Q)\in\mathfrak H_{p,r}$ such that $F_Q=F^h$ (see (\ref{Fh1})).
Moreover, the mapping $\Upsilon:\mathfrak Q_p\to\mathfrak H_{p,r}$ is locally Lipschitz.
\end{theorem}
\begin{proof} Firstly, note that the mapping $\mathfrak Q_p\owns Q\mapsto F_Q\in G_{p,2r}$ is locally Lipschitz. Indeed, in view of Lemma~\ref{Kprop} one has that the mapping $\mathfrak Q_p\owns Q\mapsto K_Q\in  G_{p,2r}^+$ is locally Lipschitz. Taking into account Proposition~\ref{pr.6a}, we then easily find that the mapping $\mathfrak Q_p\owns Q\mapsto F_Q\in G_{p,2r}$ is locally Lipschitz as well.

Assume that for each $Q\in\mathfrak Q_p$ there is $h\in L_p((-1,1),\mathcal M_r)$ such that $F_Q=F^h$. Evidently, such $h$ is unique and one has $h=\eta(F),$ where $\eta:G_{p,2r}\to L_p((-1,1),\mathcal M_r)$ is a continuous linear mapping acting by the formula
$$
[\eta(F)](x):=\left\{\begin{array}{cl}F_{21}(-2x-1,1),&-1\le x\le-\frac12,\\
F_{11}(2x+1,1),&-\frac12<x\le0,\\
F_{22}(-2x+1,1),&0<x\le\frac12,\\
F_{12}(2x-1,1),&\frac12<x\le1,\end{array}\right.
$$
where
\begin{equation}\label{FB}
F=\begin{pmatrix}F_{11}&F_{12}\\F_{21}&F_{22}\end{pmatrix},\qquad F_{ij}\in G_{p,r}.
\end{equation}
In view of (\ref{F}), note that $\mathscr F_Q\in\widetilde{\mathscr G}_{p,2r}$. Since
$\mathscr F_Q=\mathscr F^h$, it then follows from Lemma~\ref{FhLemma} that $h\in\mathfrak H_{p,r}$.
Moreover, since the mapping $\mathfrak Q_p\owns Q\mapsto F_Q\in G_{p,2r}$ is locally Lipschitz, it follows that the mapping $\mathfrak Q_p\owns Q \mapsto \Upsilon(Q):=\eta(F_Q)\in\mathfrak H_{p,r}$ is locally Lipschitz as well.

Therefore, Theorem~\ref{UpsTh} will be proved if we show that for each $Q\in\mathfrak Q_p$ there is $h\in L_p((-1,1),\mathcal M_r)$ such that $F_Q=F^h$. Obviously, it suffices to prove this only for smooth functions $Q$.

So, let $Q\in\mathfrak Q_p\cap C^1([0,1],\mathcal M_{2r})$ and $F:=F_Q$.
It then follows from Proposition~\ref{appProp} that $F\in C^1(\overline\Omega_\pm,\mathcal M_{2r})$ and that
\begin{align}\label{AF}
JF'_x(x,t)+F'_t(x,t)J=0,\qquad (x,t)\in\overline\Omega_\pm,
\\\label{Fa}
F(x,0)a^*=0,\quad aF(0,x)=0,\qquad x\in(0,1).
\end{align}
If we write $F$ in the block form (\ref{FB}), we then obtain from (\ref{AF}) that
$$
\left(\tfrac\partial{\partial x}+\tfrac\partial{\partial t}\right) F_{11}=
\left(\tfrac\partial{\partial x}+\tfrac\partial{\partial t}\right) F_{22}=0,\quad
\left(\tfrac\partial{\partial x}-\tfrac\partial{\partial t}\right) F_{12}=
\left(\tfrac\partial{\partial x}-\tfrac\partial{\partial t}\right) F_{21}=0.
$$
Therefore, it follows that $F$ can be written in the form
$$
F(x,t)=\frac12\begin{pmatrix}h_1\left(\tfrac{x-t}2\right)&h_2\left(\tfrac{x+t}2\right)\\
h_3\left(-\tfrac{x+t}2\right)&h_4\left(-\tfrac{x-t}2\right)\end{pmatrix},\qquad (x,t)\in \overline\Omega_\pm,
$$
where $h_1,h_4\in C([-1/2,1/2],\mathcal M_r)$, $h_2\in C([0,1],\mathcal M_r)$ and
$h_3\in C([-1,0],\mathcal M_r)$.
Next, we find from (\ref{Fa}) that $h_1=h_4$ and that
\begin{align*}&
h_2(x)=h_1(x),\qquad x\in\left[0,\tfrac12\right],
\\&
h_3(x)=h_1(x),\qquad x\in\left[-\tfrac12,0\right].
\end{align*}
We then arrive at $F=F^h$ with $h$ given by the formula
$$
h(x):=\begin{cases}h_1(x),&x\in\left(-\tfrac12,\tfrac12\right),\\
h_2(x),&x\in\left(\tfrac12,1\right),\\
h_3(x),&x\in\left(-1,-\tfrac12\right).
\end{cases}
$$
\end{proof}

\subsection{Proof of Theorem~\ref{mainTh}}

From Lemma~\ref{dirLip} we already know that the Krein mapping $\Theta$ is locally Lipschitz. Since the mapping $\Upsilon$ from Theorem~\ref{UpsTh} is also locally Lipschitz, Theorem~\ref{mainTh} will be proved if we show that $\Upsilon=\Theta^{-1}$. Since $\mathfrak Q_p\cap C^1([0,1],\mathcal M_{2r})$ is dense everywhere in $\mathfrak Q_p$ and $\mathfrak H_{p,r}\cap C^1([-1,1],\mathcal M_r)$ is dense everywhere in $\mathfrak H_{p,r}$, it suffices to prove the equalities
\begin{equation}\label{TU}
\Theta(\Upsilon(Q))=Q, \qquad Q\in\mathfrak Q_p\cap C^1([0,1],\mathcal M_{2r}),
\end{equation}
\begin{equation}\label{UT}
\Upsilon(\Theta(h))=h, \qquad h\in\mathfrak H_{p,r}\cap C^1([-1,1],\mathcal M_r).
\end{equation}

Let us first prove (\ref{TU}). Let $Q\in\mathfrak Q_p\cap C^1([0,1],\mathcal M_{2r})$ and $h:=\Upsilon(Q)$. Since by virtue of the definition of the mapping $\Upsilon$ one has $F_Q=F^h$, in view of formula (\ref{F}) one has
$$
\mathscr I+\mathscr F^h=(\mathscr I+\mathscr K_Q)^{-1}(\mathscr I+\mathscr K_{Q^*}^*)^{-1}.
$$
From the other hand, we obtain from Remark~\ref{rem2} that
$$
F^h(x,t)+L_h(x,t)+\int_0^x L_h(s,t)F^h(x,s)\,{\mathrm d} s=0,\qquad (x,t)\in\overline\Omega_+,
$$
where $L_h\in G_{p,2r}^+$ is of (\ref{LH}).
In view of Remark~\ref{rem1} we then find that $K_Q=L_h$. By virtue of formulas (\ref{LH}) and (\ref{JK}), it then holds
$$
[\Theta(h)](x)=K_Q(x,x)J-JK_Q(x,x)=Q(x),\qquad x\in(0,1),
$$
as desired.

It thus only remains to prove (\ref{UT}). Let $h\in\mathfrak H_{p,r}\cap C^1([-1,1],\mathcal M_r)$ and $Q:=\Theta(h)$. Then (\ref{UT}) will be proved if we show that
\begin{equation}\label{FQh}
F_Q=F^h.
\end{equation}
In turn, since
$$
\mathscr I+\mathscr F_Q=(\mathscr I+\mathscr K_Q)^{-1}(\mathscr I+\mathscr K_{Q^*}^*)^{-1},
$$
we find from Remarks~\ref{rem2} and \ref{rem1} that (\ref{FQh}) will be proved if we show that $K_Q=L_h$ with $L_h$ of (\ref{LH}). For this purpose, it suffices to verify that the function
\begin{equation}\label{phiL}
\varphi(x,\lambda):=\varphi_0(x,\lambda)+\int_0^x L_h(x,t)\varphi_0(t,\lambda)\,{\mathrm d} t,\qquad
x\in[0,1],\quad \lambda\in\mathbb C,
\end{equation}
where
$\varphi_0(x,\lambda):=({\mathrm e}^{{\mathrm i}\lambda x},\,{\mathrm e}^{-{\mathrm i}\lambda x})^\top$,
solves the Cauchy problem
\begin{equation}\label{verProbl}
J\frac{\mathrm d}{\mathrm d x}\varphi+Q\varphi=\lambda\varphi,\qquad
\varphi(0,\lambda)=\begin{pmatrix}I\\I\end{pmatrix}.
\end{equation}
The verification of this claim repeats the proof of Theorem~3.1 in \cite{dirac1}.

Indeed, let $H$ and $R_H$ be as in (\ref{H}) and (\ref{H2}), respectively. In view of Remark~\ref{rem2}, it then holds
\begin{equation}\label{RH2}
R_H(x,t)+H(x-t)+\int_0^x R_H(x,s)H(s-t)\,{\mathrm d} s=0,\qquad (x,t)\in\overline\Omega_+.
\end{equation}
Moreover, it follows from \cite[Lemma~3.4]{sturm} that in the case of the smooth $h$ as chosen one has $R_H\in C^1(\overline{\Omega}_+,\mathcal M_{2r})$.

Taking into account formulas (\ref{LH}) and $B\varphi_0(x,\lambda)=\varphi_0(-x,\lambda)$, we can rewrite (\ref{phiL}) in the form
$$
\varphi(x,\lambda)=\varphi_0(x,\lambda)+\int_0^x R_H(x,x-t)\varphi_0(x-2t,\lambda)\,{\mathrm d} t.
$$
From this equality, taking into account that
$J\tfrac{\mathrm d}{\mathrm d x}\varphi_0(x,\lambda)-\lambda\varphi_0(x,\lambda)=0$, we find that
\begin{equation}\label{aux1}
\begin{split}
J\frac{\mathrm d}{\mathrm d x}\varphi(x,\lambda)+Q(x)\varphi(x,\lambda)-\lambda\varphi(x,\lambda)
=\left\{JR_H(x,0)B\varphi_0(x,\lambda)+Q(x)\varphi_0(x,\lambda)\right\}
\\
+\int_0^x\left\{J\frac\partial{\partial x}R_H(x,x-t)
+Q(x)R_H(x,x-t)\right\}\varphi_0(x-2t,\lambda)\,{\mathrm d} t.
\end{split}
\end{equation}
Since $Q(x)=-JR_H(x,0)B$, (\ref{aux1}) is reduced to
\begin{align*}&
J\frac{\mathrm d}{\mathrm d x}\varphi(x,\lambda)+Q(x)\varphi(x,\lambda)-\lambda\varphi(x,\lambda)
\\&\qquad\qquad
=J\int_0^x\left\{\frac\partial{\partial x}R_H(x,x-t)
-R_H(x,0)BR_H(x,t)B\right\}\varphi_0(x-2t,\lambda)\,{\mathrm d} t.
\end{align*}
Therefore, (\ref{verProbl}) will be verified if we show that
\begin{equation}\label{verR}
\frac\partial{\partial x}R_H(x,x-t)-R_H(x,0)BR_H(x,t)B=0,\qquad (x,t)\in\overline\Omega_+.
\end{equation}

Let us prove (\ref{verR}). For this purpose, we obtain from (\ref{RH2}) that
$$
R_H(x,x-t)+H(t)+\int_0^x R_H(x,x-s)H(t-s)\,{\mathrm d} s=0,\qquad (x,t)\in\overline\Omega_+.
$$
Differentiating this expression in $x$ we find that
\begin{equation}\label{dRH}
\begin{split}&
\frac\partial{\partial x}R_H(x,x-t)+R_H(x,0)H(t-x)
\\&\qquad\qquad
+\int_0^x \frac\partial{\partial x}[R_H(x,x-s)]H(t-s)\,{\mathrm d} s=0,\qquad (x,t)\in\overline\Omega_+.
\end{split}
\end{equation}
Multiplying now (\ref{RH2}) by $R_H(x,0)B$ from the left and by $B$ from the right and subtracting it from (\ref{dRH}), in view also of the relation $H(x)B=BH(-x)$, we find that the function
$$
X(x,t):=\frac\partial{\partial x}R_H(x,x-t)-R_H(x,0)BR_H(x,t)B,\qquad (x,t)\in\overline\Omega_+,
$$
solves the equation
$$
X(x,t)+\int_0^x X(x,s)H(t-s)\,{\mathrm d} s=0,\qquad (x,t)\in\overline\Omega_+.
$$
Since $R_H\in C^1(\overline{\Omega}_+,\mathcal M_{2r})$, one has
$X\in C(\overline{\Omega}_+,\mathcal M_{2r})$ and thus by virtue of Proposition~\ref{hCor} we find that $X(x,t)=0$, $(x,t)\in\overline\Omega_+$. Therefore, (\ref{verR}) follows and the proof is complete.

\appendix

\section{Proof of equalities (\ref{AF}) and (\ref{Fa})}

The aim of this appendix is to prove equalities (\ref{AF}) and (\ref{Fa}) which were used in the proof of Theorem \ref{UpsTh}. The proof is technical and goes back to the well known fact that kernels of transformation operators satisfy some differential equations.

Let $\mathcal A$ be the differential operator acting on functions $X:(x,t)\mapsto\mathcal M_{2r}$ from the class $C^1(\overline{\Omega}_\pm,\mathcal M_{2r})$ by the formula
\begin{equation}\label{A}
\mathcal AX:=JX'_x+X'_tJ,
\end{equation}
where $X'_x$ and $X'_t$ denote the derivatives in variables $x$ and $t$, respectively. We shall prove the following proposition:

\begin{proposition}\label{appProp}
Let $Q\in\mathfrak Q_p\cap C^1([0,1],\mathcal M_{2r})$ and $F:=F_Q$. Then $F\in C^1(\overline{\Omega}_\pm,\mathcal M_{2r})$ and
\begin{align}\label{AF1}
(\mathcal AF)(x,t)=0,\qquad (x,t)\in \overline{\Omega}_\pm,
\\\label{Fa1}
F(x,0)a^*=0,\quad aF(0,x)=0,\qquad x\in(0,1),
\end{align}
where $a:=\begin{pmatrix}I,&-I\end{pmatrix}$.
\end{proposition}

\noindent The proof of Proposition~\ref{appProp} will be based on two auxiliary lemmas:

\begin{lemma}\label{lemma1}
Let $Q\in\mathfrak Q_p\cap C^1([0,1],\mathcal M_{2r})$ and $K:=K_Q$. Then
$K\in C^1(\overline{\Omega}_+,\mathcal M_{2r})$ and
\begin{align}\label{AK}
(\mathcal AK)(x,t)=-Q(x)K(x,t),\qquad (x,t)\in\Omega_+ .
\\\label{JK}
(KJ-JK)(x,x)=Q(x),\quad K(x,0)a^*=0,\qquad  x\in(0,1).
\end{align}
\end{lemma}

\begin{proof}
Let $Q$ and $K$ be as in the statement of the lemma. Recall (see (\ref{K})) that
\begin{equation}\label{K1}
K(x,t)=\tfrac12\left\{P^+\left(x,\tfrac{x-t}2\right)+P^+\left(x,\tfrac{x+t}2\right)B
+P^-\left(x,\tfrac{x-t}2\right)B+P^-\left(x,\tfrac{x+t}2\right)\right\},
\end{equation}
where $P^\pm$ are from (\ref{Y}) and $B$ is from (\ref{LH}). It follows from the results of \cite{cauchy} that if $Q\in\mathfrak Q_p\cap C^1([0,1],\mathcal M_{2r})$, then $P^\pm\in C^1(\overline{\Omega}_+,\mathcal M_{2r})$ and, moreover,
\begin{align}\label{auxPplus}
&P^+(x,t)=\int_t^x JQ(s)P^-(s,s-t)\,{\mathrm d} s,
\\\label{auxPminus}
&P^-(x,t)=\int_t^x JQ(s)P^+(s,s-t)\,{\mathrm d} s+JQ(t),
\\\label{auxPJ}
&P^+(x,t)J=JP^+(x,t),\qquad P^-(x,t)J=-JP^-(x,t).
\end{align}
Using now (\ref{K1})~--~(\ref{auxPJ}) and the equalities
\begin{equation}\label{Jrel}
J^2=-I_{2r}, \qquad JB=-BJ, \qquad JQ(x)=-Q(x)J,
\end{equation}
by virtue of straightforward (but quite extensive) verification we then arrive at (\ref{AK}).

Now let us prove (\ref{JK}). It follows from (\ref{auxPplus}) and (\ref{auxPminus}) that
$P^+(x,x)=0$ and $P^-(x,x)=JQ(x)$.
Therefore, in view of (\ref{K1}), we find that
$$
K(x,x)=\tfrac12\left\{P^+(x,0)+P^-(x,0)B+JQ(x)\right\}.
$$
Taking into account (\ref{auxPJ}) and  (\ref{Jrel}) we then obtain that
$$
K(x,x)J- JK(x,x)=Q(x).
$$
Since $(I_{2r}+B)a^*=0$, in view of formula (\ref{K1}) we then arrive at
$$
K(x,0)a^*=\tfrac12\left\{P^+\left(x,\tfrac{x}2\right)+P^-\left(x,\tfrac{x}2\right)\right\}
(I_{2r}+B)a^*=  0
$$
and thus (\ref{JK}) is proved.
\end{proof}

\begin{lemma}\label{lemma2}
Let $Q\in\mathfrak Q_p\cap C^1([0,1],\mathcal M_{2r})$ and $L:=L_Q$. Then
$L\in C^1(\overline{\Omega}_+,\mathcal M_{2r})$ and
\begin{align}\label{AL}
(\mathcal A L)(x,t)=L(x,t)Q(t),\qquad (x,t)\in \overline{\Omega}_+,
\\\label{JL}
(JL-LJ)(x,x)=Q(x),\quad L(x,0)a^*=0,\qquad x\in [0,1].
\end{align}
\end{lemma}

\begin{proof}
Let $Q\in\mathfrak Q_p\cap C^1([0,1],\mathcal M_{2r})$, $K:=K_Q$ and $L:=L_Q$.  In view of (\ref{L}), it follows that
$
(\mathscr I+\mathscr K_Q)(\mathscr I+\mathscr L_Q)=(\mathscr I+\mathscr L_Q)(\mathscr I+\mathscr K_Q)=\mathscr I
$\,\,
and thus for $(x,t)\in\overline{\Omega}_+$ it holds
\begin{equation}\label{auxKL1}
\begin{split}
 K(x,t)+L(x,t)+\int_t^x K(x,s)L(s,t)\,{\mathrm d} s  & =0,  \\
 K(x,t)+L(x,t)+\int_t^x L(x,s)K(s,t)\,{\mathrm d} s  & =0.   \\                                                                                                                         \end{split}
\end{equation}
In view also of (\ref{JK}), these equalities easily lead us to (\ref{JL}).

To prove (\ref{AL}), set
\begin{equation}\label{X}
S(x,t):=\int_t^x K(x,s)L(s,t)\,{\mathrm d} s,\qquad (x,t)\in\overline{\Omega}_+.
\end{equation}
Taking into account (\ref{AK}), it can be verified that
\begin{multline}\label{eq1}
(\mathcal A S)(x,t) +Q(x)S(x,t) = \\JK(x,x)L(x,t) -K(x,t)L(t,t)J -\int_t^x  K'_s(x,s)JL(s,t)\,{\mathrm d} s +
\int_t^x  K(x,s) L'_t(s,t)J\,{\mathrm d} s.
\end{multline}
Integrating by parts then leads to
$$
\int_t^x  K'_s(x,s)JL(s,t)\,{\mathrm d} s =K(x,x)JL(x,t) -K(x,t)JL(t,t) -
\int_t^x  K(x,s)J L'_s(s,t)\,{\mathrm d} s .
$$
Therefore, taking into account (\ref{JK}) and (\ref{JL}), we can rewrite (\ref{eq1}) in the form
\begin{multline}\label{eq1a}
(\mathcal A S)(x,t) +Q(x)S(x,t) = \\-Q(x)L(x,t) +K(x,t)Q(t) +\int_t^x  K(x,s)(\mathcal A L)(s,t)\,{\mathrm d} s .
\end{multline}
Now let
\begin{equation*}\label{X}
X(x,t):=(\mathcal AL)(x,t)-L(x,t)Q(t),\qquad (x,t)\in\overline{\Omega}_+,
\end{equation*}
and $X(x,t):=0$, $(x,t)\in\Omega_-$. Since $S(x,t)=-K(x,t)-L(x,t)$, from (\ref{eq1a}) and (\ref{AK}) we then find that
$$
X(x,t)+\int_0^x K(x,s)X(s,t)\,{\mathrm d} s=0.
$$
Since the operator $\mathscr I+\mathscr K$ is invertible in $\mathbb H$, we then obtain that $X(x,t)=0$ for all $(x,t)\in\overline{\Omega}_+$ which proves (\ref{AL}).
\end{proof}

Now we are ready to prove Proposition~\ref{appProp}:

\begin{proofpropa1}
Let $Q\in\mathfrak Q_p\cap C^1([0,1],\mathcal M_{2r})$,\,\, $L:=L_Q$, $L_*:=L_{Q^*}$ and $F:=F_Q$. It then follows from (\ref{L}) and (\ref{F}) that
$$
F(x,t)=L(x,t)+L_*(t,x)^*+\int_0^1 L(x,s)L_*(t,s)^*\,{\mathrm d} s,\qquad x,t\in[0,1].
$$
Since $L(x,t)=L_*(x,t)=0$ as $x<t$, we then obtain that
\begin{equation}\label{AF2}
\begin{split}
 F(x,t)=  & L(x,t)+\int_0^t L(x,s)L_*(t,s)^*\,{\mathrm d} s,\qquad  (x,t)\in {\Omega}_+,   \\
 F(x,t)=  & L_*(t,x)^*+\int_0^x L(x,s)L_*(t,s)^*\,{\mathrm d} s,\qquad (x,t)\in {\Omega}_-,                                                                                                                        \end{split}
\end{equation}
which immediately implies (\ref{Fa1}). Furthermore, it follows from (\ref{AF2}) and Lemma~\ref{lemma2} that $F\in C^1(\overline{\Omega}_\pm,\mathcal M_{2r})$.

To prove also (\ref{AF1}), take into account (\ref{F}) and observe that
$\mathscr F^*:=\mathscr F_Q^*=\mathscr F_{Q^*}$.
Therefore, it suffices to prove (\ref{AF1}) only for $(x,t)\in {\Omega}_+$. Taking into account (\ref{AL}), we obtain from the first equality in (\ref{AF2}) that
\begin{multline} \label{eq5}
(\mathcal A F)(x,t)=(\mathcal A L)(x,t)
+L(x,t)L_*(t,t)^*J -\int_0^t L'_s (x,s) JL_*(t,s)^*\,{\mathrm d} s\\
-\int_0^t L(x,s)Q(s)L_*(t,s)^*\,{\mathrm d} s +
\int_0^t L(x,s)Q(s)L_*(t,s)^*\,{\mathrm d} s +
 \int_0^t L(x,s)[(L_*(t,s))'_sJ]^*\,{\mathrm d} s.
\end{multline}
Integrating by parts leads to
\begin{multline}\label{eq6}
\int_0^t L'_s (x,s) JL_*(t,s)^*\,{\mathrm d} s=\\
 L(x,t)JL_*(t,t)^* -L(x,0)JL_*(t,0)^* -
\int_0^t L(x,s)J(L_*(t,s)^*)'_s \,{\mathrm d} s.
\end{multline}
Furthermore, in view of Lemma~\ref{lemma2} it holds
\begin{equation}\label{eq7}
JL_*(t,t)- L_*(t,t)J=Q(t)^*, \qquad t\in [0,1].                                                                                                                         \end{equation}
Taking into account (\ref{eq6}), (\ref{eq7}) and (\ref{AL}), we then obtain from (\ref{eq5}) that
\begin{equation*}\label{eq.8}
(\mathcal A F)(x,t)=
(\mathcal A L)(x,t) - L(x,t)Q(t)+L(x,0)JL_*(t,0)^*=L(x,0)JL_*(t,0)^* .
\end{equation*}
Finally, noting that $J=a^*aJ+Ja^*a$ and, in view of (\ref{JL}),
$$
L(x,0)a^*=0=L_*(x,0)a^*,\qquad x\in[0,1],
$$
we then find that $L(x,0)JL_*(t,0)^*=0$, $(x,t)\in\Omega_+$, which completes the proof of the proposition.
\end{proofpropa1}

\end{document}